\documentclass[a4paper,11pt,leqno]{article}
\usepackage{amsmath,amssymb,epsfig,psfrag,mathrsfs,dsfont}

\usepackage{hyperref}

\newtheorem {thm}{Theorem}[section]
\newtheorem {prop}[thm]{Proposition}
\newtheorem {lem}[thm]{Lemma}
\newtheorem {cor}[thm]{Corollary}
\newtheorem {rem}[thm]{Remark}
\newtheorem {example}[thm]{Example}

\newcommand{\Definition}{\bigskip\noindent{\bfseries Definition. }}

\newcommand{\Proof}[1][]{\noindent{\itshape Proof#1: }}
\newcommand{\EndProof}{~$\Diamond$}

\numberwithin{equation}{section}

\def\ba{\begin{array}}
\def\ea{\end{array}}
\def\be{\begin{equation} \label}
\def\ee{\end{equation}}
\def\bea{\begin{eqnarray*}}
\def\eea{\end{eqnarray*}}
\def\beal{\begin{eqnarray} \label}
\def\eeal{\end{eqnarray}}
\def\bit{\begin{itemize}}
\def\eit{\end{itemize}}
\def\ben{\begin{enumerate}}
\def\een{\end{enumerate}}

\def\inv{^{-1}}
\def\ti{\to\infty}
\def\#{\text{\rm card\hspace{1pt}}}
\def\1{\mathds{1}}
\def\per{\text{\rm per}}

\def\hc{\text{\rm hc}}
\def\tmp{\text{\rm tp}}
\def\pr{\text{\rm pr}}

\def\ro{\/\mathclose{[}} 

\def\RR{\mathbb{R}}
\def\NN{\mathbb{N}}
\def\ZZ{\mathbb{Z}}

\def\QQ{\mathbb{Q}}

\def\a{\alpha}
\def\b{\beta}
\def\g{\gamma}
\def\d{\delta}
\def\e{\varepsilon}
\def\l{\lambda}
\def\th{\vartheta}
\def\ph{\varphi}
\def\r{\varrho}
\def\s{\sigma}
\def\t{\tau}
\def\o{\omega}

\def\D{\Delta}
\def\L{\Lambda}
\def\Th{\Theta}
\def\O{\Omega}

\def\sA{\mathsf{A}}
\def\sB{\mathsf{B}}
\def\sC{\mathsf{C}}
\def\sD{\mathsf{D}}
\def\sS{\mathsf{S}}

\def\cF{\mathcal{F}}

\def\cL{\mathcal{L}}
\def\cM{\mathcal{M}}
\def\cR{\mathcal{R}}

\def\cT{\mathcal{T}}

\def\sG{\mathscr{G}}
\def\sM{\mathscr{M}}
\def\cP{\mathscr{P}}
\def\Tau{\mathscr{T}}

\def\br{\mathbf{r}}

\begin{document}

\title{\textbf{
Variational characterisation of Gibbs measures with Delaunay triangle interaction }
} 
 
\author{David Dereudre\footnotemark[1] \ and Hans-Otto Georgii\footnotemark[2]}
\date{\today}
\maketitle

\begin{abstract} 
This paper deals with stationary Gibbsian point processes on the plane 
with an interaction that depends on the tiles of the Delaunay triangulation of points 
via a bounded triangle potential. It is shown that the class of these Gibbs 
processes includes all minimisers of the associated free energy density and is therefore nonempty. Conversely, each such Gibbs process minimises the free energy density,
provided the potential satisfies a weak long-range assumption.\\

{\bf Running head}. {Delaunay-Gibbs measures}

\end{abstract}

\footnotetext[1]{LAMAV, Universit\'e de Valenciennes et du Hainaut-Cambr\'esis, Le Mont Houy
 59313 Valenciennes Cedex 09, France. 
E-mail: david.dereudre@univ-valenciennes.fr}

\footnotetext[2]{Mathematisches Institut der Universit\"{a}t
         M\"{u}nchen, Theresienstra{\ss}e 39, 80333 M\"{u}nchen, Germany.
         E-mail: georgii@math.lmu.de}

\section{Introduction}

It is well-known that stationary renewal processes with a reasonable spacing distribution
can be characterised as Gibbs processes for an interaction between nearest-neighbour
pairs of points \cite[Section 6]{HS}.
Here we consider an analogue in two dimensions, viz. Gibbsian point processes on 
$\mathbb{R}^2$ with an interaction depending on nearest-neighbour triples
of points, where the nearest-neighbour triples are defined in terms of the Delaunay triangulation. Recall that the Delaunay triangulation is dual to the Voronoi tessellation,
in the sense that two points are connected by a Delaunay edge if and only if their Voronoi cells have a common edge. Since the Voronoi cell of a point consists of the
part of space that is closer to this point than to any other point, this means that the
Delaunay graph defines a natural nearest-neighbour structure between the points.
(Of course, the analogy with renewal processes does not reach too far because
the independence of spacings under the Palm distribution, which is characteristic of
one-dimensional renewal processes, is lost in two dimensions due to the geometric constraints.) 

There is a principal difference between the Delaunay interactions considered here
and the pair interactions that are common in Statistical Physics. Namely, suppose
a point configuration $\o$ is augmented by a new particle at $x$. In the case of 
pair interactions, $x$ is subject to some additional interaction with
the particles in $\o$, but the interaction between the particles of $\o$ is not affected by $x$. In the Delaunay case, however, the particle at $x$ not only gives rise
to some new tiles of the Delaunay triangulation, but also
destroys some other tiles that were present in the triangulation of $\o$. 
This so-called non-hereditary nature of the Delaunay triangulation 
blurs the usual distinction between attractive and repulsive interactions and
makes it difficult to use a local characterisation
of Gibbs measures in terms of their Campbell measures and Papangelou intensities.
Such a local approach to the existence of Gibbs measures for Delaunay interactions
was used in the previous work \cite{BBD99a, BBD99b, BBD08, Dereudre}
and made it necessary to impose geometric constraints on the interaction 
by removing triangles with small angles or large circumcircles.

In this paper we address the existence problem from a global point
of view, which is based on 
stationarity and thermodynamic quantities 
such as pressure and free energy density. Specifically, we show that all
minimisers of the free energy density are Gibbsian, which implies the
existence of Delaunay-Gibbs measures because the entropy density 
has compact level sets.
The converse part of this variational principle is harder and requires the
comparison of different boundary conditions in the thermodynamic limit.
It is here that the non-hereditary nature of the interaction shows up again, 
but it can be controlled with the help of stationarity and an additional condition
which is much weaker than the geometric constraints mentioned above. 
In contrast to \cite{Dereudre}, however,
we need to assume throughout that the interaction potential is bounded, and therefore 
do not cover hard-core interactions that forbid particular shapes of the tiles.
We note, however, that some ideas of the present paper can be used to
establish the existence of Delaunay-Gibbs measures in a more general setting 
that includes also the hard-core case, see \cite{DDG}. 
The extension of the variational principle to such interactions is left to future work. 
As a final comment, let us emphasise that we make repeated use of
Euler's polyhedral formula and the resulting linear complexity of the Delaunay
triangulations, and are therefore limited to two dimensions, as was already the case
in the previous papers mentioned above.

\section{Preliminaries} 

\subsection{Configurations and Delaunay triangulations}

A subset $\o$ of $\RR^2$ is called locally finite if $\#(\o\cap\D)<\infty$ for all bounded
$\D\subset\RR^2$; each such $\o$ is called a \emph{configuration}. 
We write $\O$ for the set of all configurations $\o$. 
The configuration space $\O$ is equipped with the 
$\s$-algebra $\cF$ that is generated by the counting variables $N_\D:\o\to \#(\o\cap\D)$, with $\D$ an arbitrary bounded Borel subset of $\RR^2$. 

For each Borel set $\L\subset\RR^2$ we write $\O_\L=\{\o\in\O:\o\subset\L\}$
for the set of configurations in $\L$, $\pr_\L:\o\to\o_\L:=\o\cap\L$ for the projection
from $\O$ to $\O_\L$, $\cF_\L'= \cF|\O_\L$ for the trace $\s$-algebra of $\cF$ on 
$\O_\L$, and $\cF_\L=\pr_\L\inv\cF_\L'$ for the $\s$-algebra of events in $\O$ that happen in $\L$ only.

For each configuration $\o\in\O$ we consider the Delaunay triangulation $\sD(\o)$ associated to $\o$. By definition, 
\be{eq:sD}
\sD(\o)=\big\{ \tau\subset\o: \ \#\tau=3,\, \o\cap B(\tau)=\emptyset\big\}\,,
\ee
where $B(\tau)$ is the unique open disc with $\tau\subset\partial B(\tau)$. 
$\sD(\o)$ is uniquely defined and determines a triangulation
of the convex hull of $\o$ whenever
$\o$ is in general circular position, in that no four points of $\o$ lie on a circle that contains no further points of $\o$ inside  \cite{Moller}.
If this is not the case,  one can apply some determistic rule to make the Delaunay triangulation unique. Indeed, let 
\[
\cT:=\big\{\tau\subset\RR^2:\;\#\tau=3\big\}=\big\{(x,y,z)\in(\RR^2)^3:x\prec y\prec z\big\}
\]
be the set of all triangles (or tiles) in $\RR^2$ where `$\prec$' stands for the lexicographic order in $\RR^2$. The triangles in $\cT$ can be compared by the lexicographic order
of $(\RR^2)^3$, and this in turn induces a lexicographic order on finite collections
of triangles. Now, if $n\ge4$  points of $\o$ lie on a circle with no points inside then the associated Delaunay cell is a convex polygon having these $n$ points as vertices. 
To define a unique triangulation of this polygon one can then simply take
the smallest among all possible triangulations. 
Conflicts with other possible polygons cannot arise because the
tessellations inside and outside a fixed convex polygon do not depend on each other. 

Let us note that the prescription $\o\to\sD(\o)$ is a mapping from $\O$ to the set 
$\O(\cT)$ of all locally finite subsets of $\cT$. If $\O(\cT)$ is equipped with the 
$\s$-algebra $\cF(\cT)$ that is defined in analogy to $\cF$, one can easily 
check that this mapping is measurable.

Next we assign to each tile $\tau\in\cT$ a centre and a radius. Specifically,
for every $\tau\in\cT$ we write $c(\tau)$ for the centre and $\r(\tau)$ for the radius of 
the circumscribed disc $B(\tau)$. 
The centres allow us to consider $\sD(\o)$ as a germ-grain system, i.e., as
a marked point configuration of germs in $\RR^2$ and marks in the space
\[
\cT_0= \{ \tau\in\cT: \;c(\tau)=0\}
\]
of centred tiles, by considering the mapping
\be{eq:D}
D: \o \to \big\{(c(\tau),\tau-c(\tau)): \tau\in\sD(\o)\big\}
\ee
from $\O$ to the point configurations on $\RR^2\times\cT_0$. Here we write
$\tau-c(\tau):=\{y-c(\tau):y\in\tau\}$ for the shifted tile.

A crucial fact we need in the following is the linear complexity of 
Delaunay triangulations, which is expressed in the following
lemma. This result follows directly from Euler's polyhedral formula, and is the 
main reason why we need to confine ourselves to two spatial dimensions;
see \cite{Proofs}, Chapter 11, and \cite{Moller}, Remark 2.1.4.

\begin{lem}\label{Euler}
For a simple planar graph on $n$ vertices, the number of edges is at most
$3n-6$, and the number of inner faces is at most $2n-5$. In particular, 
every triangulation with $n$ nodes consists of $2n-2-\partial$ triangles, 
where $\partial$ is the number of nodes (or: number of edges) along 
the outer boundary.
\end{lem} 

\subsection{Stationary point processes and their tile distribution}

Let $\cP_\Th$ be the set of all probability measures $P$ on $(\O,\cF)$ that satisfy
the following two properties:
\bit
\item[(S)] $P$ is \emph{stationary}, that is, $P$ is invariant under the shift group 
$\Th=(\th_x)_{x\in\RR^2}$ on $\O$, which is defined by $\th_x:\o\to\o-x:=\{y-x:y\in\o\}$.
\item[(I)] $P$ has a \emph{finite intensity} $z(P)=|\D|\inv\int N_\D\,dP<\infty$.
Here, $\D$ is any bounded Borel set in $\RR^2$ (which can be arbitrarily chosen due to stationarity), and $|\D|$ its Lebesgue measure.
\eit
Each $P\in\cP_\Th$ is called a \emph{stationary point process on $\RR^2$ with finite intensity}. For $\L\subset\RR^2$, we write $P_\L:=P\circ\pr_\L\inv $ for the projection
image of $P$ on $(\O_\L,\cF_\L')$, which can of course be identified with the restriction 
$P|\cF_\L$ of $P$ to the events in $\L$.

Every $P\in \cP_\Th$ defines a germ-grain model $\bar{P}$, namely the distribution of $P$ under the mapping $D$ defined in \eqref{eq:D}.
That is,  $\bar{P}$ is a stationary marked point process on $\RR^2$ with mark space $\cT_0$. Let  $\bar{P}^0$ be the associated Palm measure on $\cT_0\times\O$
and $\mu_P=\bar{P}^0(\,\cdot\,\times\O)$ the associated mark distribution, or \emph{centred tile distribution}, on $\cT_0$. By definition,
\begin{equation}\label{eq:Palm}
\int dx \int \mu_P(d\tau) \,f(x,\tau)
=\int P(d\o)\sum_{\tau\in\sD(\o)} f(c(\tau),\tau-c(\tau))
\end{equation}
for all nonnegative measurable functions $f$ on $\RR^2\times\cT_0$.
For each $P\in\cP_\Th$, $\mu_P$ has total mass $\|\mu_P\|=2z(P)$, as follows
from Euler's polyhedral formula; see, for example, \cite[Eq.~(3.2.11)]{Moller}
or \cite[Theorem 10.6.1(b)]{SchnWeil}. 

Let us say a measure $P\in \cP_\Th$ is \emph{tempered} if 
\be{eq:H}
\int |B(\tau)|\, \mu_P(d\tau) <\infty
\ee
We write $\cP_\Th^\tmp$ for the set of all tempered $P\in\cP_\Th$. 
Of course, \eqref{eq:H} is equivalent to the condition  
$\int \r(\tau)^2\, \mu_P(d\tau) <\infty$. Moreover, \eqref{eq:Palm} implies that
\begin{eqnarray}\label{eq:tempered2}
 &&\int |B(\tau )| \,\mu_P(d\tau) =\int dx\int \mu_P(d\tau)\, \1_{\{x\in B(\tau)\}}\\
&=&\int \sum_{\tau \in\sD(\o)} \1_{\{c(\tau)\in B(\tau-c(\tau))\}} \,P(d\o)  
=\int \#\big\{\tau\in\sD(\o): 0\in B(\tau)\big\}\, P(d\o) \,. \nonumber
\end{eqnarray}
So, $P$ is tempered if and only if the last expression is finite.
A sufficient condition for temperedness will be given in Proposition \ref{tempered}.

The most prominent members of $\cP_\Th^\tmp$ are the Poisson point processes, 
which will take the role of reference processes for the models we consider. 
Recall that the Poisson point process $\Pi^z$ with intensity $z>0$ is characterised 
by the following two properties:
\bit
\item[(P1)] For every bounded Borel set $\D$,  the counting variable $N_{\D}$ is Poisson distributed with parameter $z|\D|$. 
\item[(P2)] Conditional on $N_\D=n$, the $n$ points in $\D$ are independent with uniform distribution on $\D$, for every bounded Borel set $\D$ and each integer $n$.
\eit
The temperedness of $\Pi^z$ follows from Proposition 4.3.1 of \cite{Moller},
or Proposition \ref{tempered} below.

Another type of measures in $\cP_\Th^\tmp$ are the {stationary empirical fields}
that are defined as follows. Let $\L\subset\RR^2$ be an open   
square of side length $L$, and for $\o\in\O_\L$ let
$\o_{\L,\per}=\{x+Li:\,x\in\o,\,i\in\ZZ^2\}$
be its periodic continuation. The associated \emph{stationary empirical field}
is then given by
\be{eq:R}
R_{\L,\o}= \frac1{|\L|}\int_{\L} \d_{\th_x \o_{\L,\per}}\;dx\,.
\ee
It is clear that $R_{\L,\o}$ is stationary. In addition, it is tempered because
$2\r(\tau)\le \text{diam\,}\L$ for each triangle $\tau\in\sD(\o_{\L,\per})$.

\subsection{The topology of local convergence}

In contrast to the traditional weak topology on the set $\cP_\Th$ of stationary point processes, we exploit here a finer topology, which is such that
the intensity is a continuous function, but nonetheless the
entropy density has compact level sets. 

Let $\cL$ denote the class of all measurable functions $f:\O\to\RR$ which are
\emph{local} and \emph{tame}, in that there exists some bounded Borel set $\D\subset\RR^2$ such that $f= f\circ\pr_\D$ and $|f|\le b(1+N_\D)$ for some constant $b=b(f)<\infty$. The \emph{topology $\Tau_\cL$ of local convergence} on $\cP_\Th$ is then defined as the weak* topology induced by $\cL$, i.e., as the smallest topology for which the mappings $P\to\int f\,dP$ with $f\in\cL$ are continuous. By the definition of the intensity,
it is then clear that the mapping $P\to z(P)$ is continuous.

A further basic continuity property is the fact that the centred tile distribution $\mu_P$ depends continuously on $P$. 
Let $\cL_0$ be the class of all \emph{bounded} measurable functions on the space 
$\cT_0$  of centred tiles,  and $\Tau_0$ the associated weak* topology on the set $\cM(\cT_0)$ of all finite measures on $\cT_0$. (This is sometimes called the $\tau$-topology.) 

\begin{prop}\label{muP_cont}
Relative to the topologies $\Tau_\cL$ and $\Tau_0$ introduced above, the mapping 
$P\to\mu_P$ from $\cP_\Th$ to $\cM(\cT_0)$ is continuous.
\end{prop}
This result will be proved in Section~\ref{subsec:energy}. It takes advantage of the linear complexity of finite Delaunay triangulations, and therefore relies on the planarity of our model.

\subsection{The entropy density}

Given a point process $P\in\cP_\Th$ and a bounded Borel set $\L$ in $\RR^2$,
we let $I_\L(P,\Pi^z)$ denote the relative entropy (or Kullback-Leibler information)
of $P_\L$ relative to $\Pi^z_\L$. By the independence properties of $\Pi^z$,
these quantities are subadditive in $\L$, which implies that the limit
\be{eq:entropy}
I_z(P) =\lim_{|\L|\ti}I_\L(P;\Pi^z)/|\L|\in[0,\infty]
\ee
exists and is equal to the supremum of this sequence. For our purposes, it is
sufficient to take this limit along a fixed sequence of squares; for example,
one can take squares with vertex coordinates in $\ZZ+1/2$. 
The claim then follows from the well-known analogous result for lattice models 
\cite[Chapter 15]{GiiGibbs} by dividing $\RR^2$ into unit squares.
$I_z$ is called the (negative) \emph{entropy density with reference measure $\Pi^z$.}

We set $I=I_1$. Each $I_z$ differs from $I$ only by a constant and a multiple of the particle density. In fact, an easy computation shows that
\be{eq:I_z}
I(P) = I_z(P) + 1-z+z(P) \log z\qquad \text{for all $z>0$ and $P\in\cP_\Th$\,.}
\ee
A crucial fact we need later is the following result obtained in Lemma 5.2 of~\cite{GiiZess}.
\begin{lem}\label{entropy_cp}
In the topology~$\Tau_\cL$, each $I_z$ is lower semicontinuous with compact level sets
$\{I_z\le c\}$, $c>0$.
\end{lem}

\subsection{Triangle interactions}

This paper is concerned with point processes with a particle interaction which is
induced by the associated Delaunay triangulation. We stick here to the simplest kind
of interaction, which depends only on the triangles that occur in each configuration.
Specifically, let $\ph: \cT_0\to\RR$ be an arbitrary measurable function. 
It can be extended to a unique shift-invariant measurable function $\ph$ on $\cT$ via 
$\ph(\tau):= \ph(\tau-c(\tau))$, $\tau\in\cT$. 
Such a $\ph$ will be called a \emph{triangle potential}.
We will assume throughout that $\ph$ is bounded, in that 
\begin{equation}\label{eq:ph_lb}
|\ph| \le c_\ph
\end{equation}
for some constant $c_\ph<\infty$.  
In Theorem \ref{Gtempered} we will need the following additional condition
to prove the temperedness of Gibbs measures.
Let us say that a triangle potential $\ph$ is \emph{eventually increasing} if there exist a constant
$r_\ph<\infty$ and a measurable nondecreasing function $\psi:[r_\ph,\infty\ro\to\RR$
such that $ \ph(\tau) = \psi(\r(\tau))$ when \mbox{$\r(\tau)\ge r_\ph$.}
This condition is clearly satisfied when $\ph$ is constant for all triangles $\tau$ with
sufficiently large radius $\r(\tau)$.
  
\begin{example}\label{examples}\rm
Here are some examples of triangle potentials.
For each triangle $\tau\in\cT$ let $b(\tau)=\frac 13\sum_{x\in\tau}x$ be the barycentre
and $A(\tau)$ the area of $\tau$. Examples of bounded (and scale invariant) 
interactions that favour equilateral Delaunay triangles are
\[
\ph_1(\tau) =\b\, |c(\tau)-b(\tau)|/\r(\tau)  \quad\text{ or }\quad \ph_2(\tau) = - \b\, A(\tau)/\r(\tau)^2
\]
with $\b>0$. Of course, many variants are possible; e.g., one can replace the 
barycentre by the centre of the inscribed circle. 
By way of contrast, to penalise regular configurations one can replace the $\ph_i$'s 
by their negative.

The triangle potentials $\ph_i$ above are not eventually increasing. 
But each triangle potential $\ph$ can be modified to exhibit this property by setting
\[
\tilde\ph(\tau) =\left\{ \ba{cl}\ph(\tau) & \text{ if }  \r(\tau)<r\,,\\
K &\text{ otherwise} \ea\right.
\]
with $r>0$ and $K$ a suitable constant. When $K$ is large, one has the additional effect of favouring small circumcircles.
\end{example}

\begin{rem}\label{rem:edges}\rm 
The type of interaction introduced above is the simplest possible 
that is adapted to the Delaunay structure. In particular, we avoid here any explicit interaction $\psi$ along the Delaunay edges. This has two reasons: 
First, we might add a term of the form
$\frac12 \sum_{e\subset\tau:\,\# e=2}\ph_\text{edge}(e)$ to the triangle interaction $\ph$.
Such a term would take account for an edge interaction
$\ph_\text{edge}$ whenever $\sD(\o)$ is a triangulation of the full plane.
Secondly, we often need to control the interaction over large distances;
the condition of $\ph$ being eventually increasing is tailored for this
purpose. It is then essential to define the
range in terms of triangles rather than edges. Namely, 
if a configuration $\o$ is augmented by a particle at a large distance from $\o$,
the circumcircles of all destroyed triangles must be large, but their edges can be arbitrarily short. So, a large-circumcircle assumption on the triangle potential allows
to control this effect, but a long-edge asumption on an edge potential would be useless.
\end{rem}

\section{Results}

Let $\ph$ be a fixed triangle potential. \emph{We assume throughout that $\ph$ is bounded}, see \eqref{eq:ph_lb},
but do not require in general that $\ph$ is eventually increasing. For each bounded Borel set $\L\subset\RR^2$ and each configuration $\o\in\O$,
the associated \emph{Hamiltonian in $\L$ with boundary condition $\o$} is then
defined for arbitrary $\zeta\in\O$ by
\be{eq:cH_o}
H_{\L,\o}(\zeta)= \sum_{\tau\in\sD(\zeta_\L\cup\o_{\L^c})
:\,B(\tau)\cap\L\ne\emptyset}\ph(\tau)\,.
\ee
It is always well-defined since the sum is finite. Note that $H_{\L,\o}(\emptyset)\ne0$ in general.  For defining the associated
Gibbs distribution we need to impose a condition on the boundary condition $\o$. 

\Definition \label{admissible} Let us say a configuration $\o\in\O$  is \emph{admissible} if
for every bounded Borel set $\L$ there exists a bounded Borel set 
$\bar\L(\o)\supset\L$ such that $B(\tau)\subset \bar\L(\o)$
whenever $\zeta\in\O$ and $\tau\in\sD(\zeta_\L\cup\o_{\L^c})$ is such that $B(\tau)\cap\L\ne\emptyset$.
We write $\O^*$ for the set of all admissible configurations. 

\bigskip\noindent
In Corollary \ref{temperedness} we will show 
that $P(\O^*)=1$ for all $P\in\cP_\Th$ with $P(\{\emptyset\})=0$. 
Suppose now that $\o\in\O^*$ and $\L$ is any bounded Borel set.
Lemma \ref{Euler} then shows that
\[
H_{\L,\o}\ge -2c_\ph N_\L -2c_\ph N_{\bar\L(\o)\setminus\L}(\o)\,,
\]
where $\bar\L(\o)$ is as above. This in turn implies that for each $z>0$ the associated \emph{partition function}
\be{eq:Z_L}
Z_{\L,z,\o} = \int e^{-H_{\L,\o}} \,d\Pi^z_{\L}
\ee
is finite. We can therefore define
the \emph{Gibbs distribution with activity $z>0$} by
\be{eq:cG_o}
G_{\L,z,\o}(A) = Z_{\L,z,\o}^{-1} \int\1_A(\zeta\cup\o_{\L^c})\, e^{-H_{\L,\o}(\zeta)} \,\Pi^z_{\L}(d\zeta)\;,\quad A\in\cF.
\ee
The measure $G_{\L,z,\o}$ depends measurably on 
$\o$ and thus defines a probability kernel from $(\O^*, \cF_{\L^c})$ to $(\O,\cF)$. 

\Definition A probability measure $P$ on $(\O,\cF)$ is called a \emph{Gibbs 
point process for the Delaunay triangle potential $\ph$ and the activity $z>0$},
or a \emph{Delaunay-Gibbs measure} for short, if  $P(\O^*)=1$ and
$P = \int P(d\o)\,G_{\L,z,\o}$ for all bounded Borel sets $\L\subset\RR^2$.
We write 
$\sG_\Th(z,\ph)$ for the set of all stationary Gibbs measures  for $\ph$ and $z$, and 
$\sG_\Th^\tmp(z,\ph)$ for the set of all \emph{tempered} stationary 
Gibbs measures; recall Eq.~\eqref{eq:H}.

\bigskip
The above definition corresponds to the classical concept of a Gibbs measure,
which is based on the location of points. We note that an alternative
concept of Gibbs measure that considers the location of Delaunay triangles
has been proposed and used by Zessin \cite{Zess02} and Dereudre \cite{Dereudre}.

Intuitively, the interaction of a configuration in $\L$ with its boundary
condition $\o$ reaches not farther
than the set $\bar\L(\o)$ above, which guarantees some kind of quasilocality.
So one can expect that a limit of suitable Gibbs distributions  $G_{\L,z,\o}$ as 
$\L\uparrow\RR^2$ should be Gibbsian, and the existence problem reduces to the question of whether such limits exist. Our approach here is to take the necessary compactness property from Lemma \ref{entropy_cp}, the compactness of level sets of the entropy density.
In fact, we even go one step further and show that the stationary Gibbs measures
are the minimisers of the free energy density. Since such minimisers exist
by the compactness of level sets, this solves in particular the existence problem.
The free energy density is defined as follows; recall the definition of the centred tile
distribution $\mu_P$ before \eqref{eq:Palm}.

\Definition The \emph{energy density} of a stationary point process $P\in\cP_\Th$ is defined by
\[
\Phi(P)=\int \ph\,d\mu_P = 
|\D|\inv\int P(d\o)\sum_{\tau\in\sD(\o):\,c(\tau)\in\D} \ph(\tau)\,,
\]
where $\D$ is an arbitrary bounded Borel set. The \emph{free energy density} of $P$
relative to $\Pi^z$ is given by $I_z(P)+\Phi(P)$.

\bigskip
The definition of $\Phi$ is justified by Proposition \ref{prop:Phi} below.
Here are some crucial facts on the free energy density, which will be proved
in Subsection \ref{subsec:energy}.

\begin{prop}\label{I+Phi_lsc}
Relative to the topology $\Tau_\cL$ on $\cP_\Th$, $\Phi$ is continuous,
and each \mbox{$I_z+\Phi $} is lower semicontinuous with compact level sets.
In particular, the set $\sM_\Th(z,\ph)$ of all minimisers of $I_z+\Phi $
is a non-empty convex compact set, and in fact a face of the simplex $\cP_\Th$.
\end{prop} 

Next we observe that the elements of $\sM_\Th(z,\ph)$ are nontrivial,
in that the empty configuration $\emptyset\in\O$ has zero probability; this result will
also be proved in Subsection~\ref{subsec:energy}.

\begin{prop}\label{prop:non-degenerate}
For all $z>0$ we have $\d_{\emptyset}\notin\sM_\Th(z,\ph)$, and thus 
$P(\{\emptyset\})=0$ for all $P\in\sM_\Th(z,\ph)$.
\end{prop}

Our main result is the following variational characterisation of Gibbs measures.

\begin{thm}\label{varprinc}
Let $\ph$ be a bounded triangle potential and let $z>0$. Then
every minimiser of the free energy density is a stationary Gibbs measure. That is, the identity $\sM_\Th(z,\ph)\subset\sG_\Th(z,\ph)$ holds. In particular, Gibbs measures exist. Conversely, every tempered stationary Gibbs measure is a minimiser of the free energy density, which means that $\sG^\tmp_\Th(z,\ph)\subset \sM_\Th(z,\ph)$.
\end{thm}

The proof will be given in Subsections \ref{subsec:var_princ1} and
\ref{subsec:var_princ2}. Theorem \ref{varprinc} raises the problem of whether all
stationary Gibbs measures are tempered. 
It is natural to expect that $\sG_\Th(z,\ph)=\sG^\tmp_\Th(z,\ph)$, but we did not succeed to prove this
in general. In fact, we even do not know whether $\sG^\tmp_\Th(z,\ph)$ is always 
non-empty. 
But we can offer the following sufficient condition,
which will be proved in Subsection \ref{Gibbstempered}.

\begin{thm}\label{Gtempered}
Suppose $\ph$ is eventually increasing and let $z>0$. Then every stationary Gibbs measure is tempered, so that $\sG^\tmp_\Th(z,\ph)=\sG_\Th(z,\ph)$.
\end{thm}
  
 Combining Theorems \ref{varprinc} and \ref{Gtempered} we arrive at the following result.

\begin{cor}\label{varprincb}
Suppose $\ph$ is bounded and eventually increasing, and let $z>0$. Then
the minimisers of the free energy density are precisely the stationary tempered Gibbs measures. That is,  $\sM_\Th(z,\ph)=\sG_\Th(z,\ph)=\sG_\Th^\tmp(z,\ph)$ for all $z>0$. 
\end{cor}

 The proof of Theorem \ref{varprinc} is based on an analysis
of the mean energy and the pressure in the infinite volume limit when
$\L\uparrow\RR^2$. For simplicity, we take this limit through a fixed reference
sequence, namely the sequence
\be{eq:L_n}\textstyle
\L_n = \/\big]-n-\frac12, n+\frac12\/\big[^{\;2}
\ee
of open centred squares.
We shall often write $n$ when we refer to $\L_n$. That is, we set $\o_n=\o_{\L_n}$,
$P_n=P_{\L_n}$, $R_{n,\o}=R_{\L_n,\o}$, $H_{n,\o}=H_{\L_n,\o}$, and so on.
We also write $v_n=|\L_n|=
(2n+1)^2$ for the Lebesgue measure of $\L_n$. Our first result justifies the above
definition of $\Phi(P)$. Besides the Hamiltonian \eqref{eq:cH_o} with configurational
boundary condition $\o$, we will also
consider the \emph{Hamiltonian with periodic boundary condition}, namely
\be{eq:H_per}
H_{n,\per}(\o):=v_n\Phi(R_{n,\o})=\sum_{\tau\in\sD(\o_{n,\per}): \,c(\tau)\in\L_n}\ph(\tau)\;.
\ee
By definition, we have $H_{n,\per}(\emptyset)=0$. 
Applying Lemma~\ref{Euler} and using \eqref{eq:ph_lb}, we see that  $|H_{n,\per}|\le v_nc_\ph\,2\,z(R_n)= 2c_\ph\,N_n$. The following result will be proved in Subsection \ref{subsec:energy}.

\begin{prop}\label{prop:Phi}
For every $P\in \cP_\Th$ we have
\[
 \lim_{n\ti} v_n\inv\int H_{n,\per}\,dP= \Phi(P)\,.
\]
Moreover, if $P$ is tempered then
\[
\lim_{n\ti}v_n\inv\int  H_{n,\o}(\o)\,P(d\o)=\Phi(P)\,.
\]
\end{prop}

Finally we turn to the pressure. Let
\[
Z_{n,z,\per} = \int e^{-H_{n,\per}} \,d\Pi^z_n
\]
be the partition function in $\L_n$ with periodic boundary condition.

\begin{prop}\label{prop:pressure} 
For each $z>0$,  the pressure
\[
p(z,\ph):=\lim_{n\ti} v_n^{-1}\log Z_{n,z,\per}
\]
exists and satisfies
\be{eq:pressure}
p(z,\ph) =  - \min_{P\in\cP_\Th}\big [I_z(P)+\Phi(P)\big]\,.
\ee
\end{prop}
\Proof This is a direct consequence of Theorem 3.1 of \cite{GiiZess}
because $\Phi$ is continuous by Proposition \ref{I+Phi_lsc}.\EndProof

\bigskip A counterpart for the partition functions with configurational boundary
conditions follows later in Proposition \ref{lowerbound}.
Let us conclude with some remarks on extensions and further results.

\begin{rem}Large deviations. \rm The following large deviation principle is valid. 
For every measurable $A\subset\cP_\Th$,
\[
\limsup_{n\ti} v_n\inv\log G_{n,z,\per}(R_n\in A) \le - \inf I_{z,\ph}(\text{\rm cl\,}A)
\]
and
\[
\liminf_{n\ti} v_n\inv\log G_{n,z,\per}(R_n\in A) \ge - \inf I_{z,\ph}(\text{\rm int\,}A)\,,
\]
where $G_{n,z,\per}=Z_{n,z,\per}\inv e^{-H_{n,\per}}\Pi^z_n$ is the Gibbs distribution in $\L_n$ with periodic boundary condition, $I_{z,\ph}=I_z+\Phi +p(z,\ph)$ is the excess 
free energy density,
and the closure cl and the interior int are taken in the topology $\Tau_\cL$. 
Since $\ph$ is bounded so that $\Phi$  is continuous, this is a direct consequence of 
Theorem 3.1 of \cite{GiiZess}.
\end{rem}

\begin{rem}\label{rem:marks}Marked particles. \rm
Our results can be extended to the case of point particles with marks, 
that is, with internal degrees of freedom. Let $E$ be any separable metric space,
which is equipped with its Borel $\s$-algebra and a reference measure $\nu$, and 
$\overline{\O}$ the set of all pairs $\bar\o=(\o,\s_\o)$ with $\o\in\O$ and $\s_\o\in E^\o$.
In place of the reference Poisson point process $\Pi^z$, one takes
the Poisson point process $\overline\Pi^z$ on $\overline{\O}$ with intensity measure 
$z\l\otimes\nu$, where $\l$ is Lebesgue measure on $\RR^2$.
For $\bar\o\in\overline{\O}$ let 
\[
\sD(\bar\o)=\big\{\bar\tau=(\tau,\s_\tau): \tau\in\sD(\o), \ \s_\tau=\s_\o\,|_{\/\tau}\big\}\,.
\]
Of course, the centre, radius and circumscribed disc of a marked triangle $\bar\tau$ 
are still defined in terms of the underlying $\tau$.
In the germ-grain representation, $\cT_0$ is replaced by the set $\overline{\cT_0}$ of all 
centred $\bar\tau$. The tile distribution $\mu_{\bar P}$ of a stationary point process
$\bar P$ on $\overline{\O}$ is a finite measure on $\overline{\cT_0}$ and is defined
by placing bars in~\eqref{eq:Palm}. A triangle potential is a bounded function $\ph$ on 
$\overline{\cT_0}$. Such a $\ph$ is eventually increasing if $ \ph(\bar\tau) = \psi(\r(\tau))$ for some nondecreasing $\psi$ when $\r(\tau)$ is large enough.
It is then easily seen that all our arguments carry over to this setting without change. 
\end{rem}

\begin{rem}\label{rem:hc}Particles with hard core. \rm
There is some interest in the case when the particles are required to have at
least some distance $r_0>0$. This is expressed by adding to the Hamiltonian
\eqref{eq:cH_o} a hard-core pair interaction term
$H_{\L,\o}^\hc(\zeta)$ which is equal to $\infty$ if $|x-y|\leq r_0$ for a pair
$\{x,y\}\subset\zeta_\L\cup\o_{\L^c}$ with $\{x,y\}\cap\L\ne\emptyset$, and zero otherwise. Equivalently, one can replace
the configuration space $\O$ by the space
\[\O^\hc=\big\{\o\in\O: |x-y|> r_0 \text{ for any two distinct }x,y\in\o\big\}\]
of all hard-core configurations. The free energy functional on $\cP_\Th$ then takes
the form $F_z^\hc:=I_z+\Phi+\Phi^\hc$, where
\[
\Phi^\hc(P)= \infty \ P^0\big(\o: 0<|x|\le r_0 \text{ for some }x\in\o\big)
= \infty \ P(\O\setminus \O^\hc)
\]
for $P\in\cP_\Th$ with Palm measure $P^0$; here we use the convention $\infty\, 0 =0$.
We claim that our results can also be adapted to this setting. In particular, 
the minimisers of $F_z^\hc$ are Gibbsian for $z$ and the combined triangle and hard-core pair interaction, and the tempered Gibbs measures for this interaction minimise $F_z^\hc$.
We will comment on the necessary modifications in Remarks \ref{rem:hc1} and \ref{rem:hc2}.
\end{rem}

Combining the extensions in the last  two remarks we can include the following
example of phase transition.

\begin{example}\label{rem:Potts}
The Delaunay-Potts hard-core model for particles
with $q\ge2$ colours. \rm In the setup of Remark \ref{rem:marks} we have
$E=\{1,\ldots,q\}$, and the triangle potential is
\[
\ph(\bar\tau)=\left\{\ba{cl}\b&\text{if 
$\r(\tau)\le  r_1$ and $\s_\tau$ is not constant,}\\
0&\text{otherwise,}\ea\right.
\]
where $\b>0$ is the inverse temperature and $r_1>0$ is an arbitrary interaction radius. 
If one adds a hard-core pair interaction with range $r_0<r_1/ \sqrt{2}$ as in Remark 
\ref{rem:hc}, this model is similar to the model considered in \cite{BBD04}.
(Instead of a triangle potential, these authors consider an edge potential along 
the Delaunay edges that do not belong to a tile $\t$ of radius $\r(\tau)>r_1$.)
Using a random cluster representation of the triangle interaction as in
\cite{Gr94} and replacing edge percolation by tile
percolation one finds that the methods of \cite{BBD04} can be adapted
to the present model. Consequently, if $z$ and $\b$ are sufficiently large,
then the simplex $\sM_\Th(z,\ph)=\sG_\Th(z,\ph)$ has at least $q$ distinct
extreme points. 
\end{example}

\section{Proofs}
\subsection{Energy and free energy}\label{subsec:energy}

We begin with the proof of Proposition \ref{muP_cont},
which states that the centred tile distribution $\mu_P$ depends continuously on $P$.
The continuity of the energy density $\Phi$ and the lower semicontinuity of the free 
energy density $I_z+\Phi$ then follow immediately.

\bigskip
\Proof[ of Proposition \ref{muP_cont}] 
Let  $(P_\a)$ be a net in $\cP_\Th$ that converges to some $P\in\cP_\Th$.
We need to show that $\int g\,d\mu_{P_\a}\to \int g\,d\mu_{P}$ for all $g\in\cL_0$.
We can assume without loss of generality that $0\le g\le 1$.

We first consider the case that $g$ has bounded support, in that $g\le\1\{\r\le r\}$
for some $r>0$. Let
$\D\subset\RR^2$ be any bounded set of Lebesgue measure $|\D|=1$. Also, let
\[
f(\o)=\sum_{\tau\in\sD(\o)}g(\tau-c(\tau))\,\1_\D(c(\tau))\,.
\]
In view of \eqref{eq:Palm}, we then have $\int g\,d\mu_Q=\int f\,dQ$ for all
$Q\in\cP_\Th$, and in particular for $Q=P_\a$ and $Q=P$.
By the bounded support property of $g$, $f$ depends only on the configuration in the
$r$-neigbourhood $\D^r:=\{y\in\RR^2:|y-x|\le r\text{ for some }x\in\D\}$ of $\D$.
That is, $f$ is measurable with respect to $\cF_{\D^r}$. Moreover, $f \le 2 N_{\D^r}$ by Lemma~\ref{Euler}, so that $f\in\cL$. In the present case, the 
result thus follows from the definition of the topology $\Tau_\cL$.

If $g$ fails to be of bounded support, we can proceed as follows.
Let $\e>0$ be given and $r>0$ be so large that  $\mu_P(\r > r)<\e$. Since
$\|\mu_{P_\a}\| = 2z(P_\a)\to  2z(P)= \|\mu_{P}\|$ and $\mu_{P_\a}(\r\le r)\to
\mu_{P}(\r\le r)$ by the argument above, we have $\mu_{P_\a}(\r> r)\to\break
\mu_{P}(\r> r)$.
We can therefore assume without loss of generality that $\mu_{P_\a}(\r> r) <\e$
for all $\a$. Using again the first part of this proof, we can thus write
\[\begin{split}
&\int g\,d\mu_{P}-\e \le \int g\,\1_{\{\r\le r\}}\,d\mu_{P} 
= \lim_\a  \int g\,\1_{\{\r\le r\}}\,d\mu_{P_\a}\\
&\le \liminf_\a  \int g \,d\mu_{P_\a} 
\le  \limsup_\a  \int g \,d\mu_{P_\a}
\\ 
& \le  \lim_\a  \int g\,\1_{\{\r\le r\}}\,d\mu_{P_\a} +\e
=  \int g\,\1_{\{\r\le r\}}\,d\mu_{P} +\e \le \int g\,d\mu_{P} +\e \,.
\end{split}\]
Since $\e$ was chosen arbitrarily, the result follows.\EndProof

\bigskip 
We now turn to the properties of the free energy density.

\bigskip
\Proof[ of Proposition \ref{I+Phi_lsc}]
As $\ph$ belongs to $\cL_0$, the continuity of $\Phi$ follows immediately from Proposition \ref{muP_cont}.
By Lemma \ref{entropy_cp}, we can also conclude that $I_z+\Phi$ is lower semicontiuous.
Moreover,  hypothesis \eqref{eq:ph_lb} implies that the level set  $\{I_z+\Phi \le c\}$ 
is contained in $\{I_z\le c+2c_\ph z(\cdot)\}$, which by  \eqref{eq:I_z} coincides 
with the compact set $\{I_{z'}\le c+z'-1\}$ for $z'=z\exp(2c_\ph)$.

Let $P=\d_\emptyset\in\cP_\Th$ be the Dirac measure at the empty configuration.
Then $\mu_P\equiv 0$ and thus $\Phi(P)=0$. On the other hand, $I_z(P)=z$. This means
that $I_z+\Phi $ is not identically equal to $+\infty$ on $\cP_\Th$ and thus, by the compactness of its level sets, attains its infimum.
To see that the minimisers form a face of $\cP_\Th$, it is sufficient to note that 
$I_z+\Phi $ is measure affine; cf.~Theorem (15.20) of \cite{GiiGibbs}.\EndProof

\bigskip
Next we show that the minimisers of the free energy are nondegenerate.
\bigskip

\Proof[ of  Proposition \ref{prop:non-degenerate}]
The second statement follows from the first because $\cM_\Th(z,\ph)$ is a face
of $\cP_\Th$. For, suppose there exists some $P\in\cM_\Th(z,\ph)$ with 
$P(\{\emptyset\})>0$. Then $\d_{\emptyset}$ appears in the ergodic 
decomposition of $P$,
which would only be possible if  $\d_{\emptyset}\in\cM_\Th(z,\ph)$.

To prove the first statement we note that
$\Phi(\Pi^u)\le c_\ph \|\mu_{\Pi^u}\| =2c_\ph u$ for all $u>0$. Therefore, if
$z>0$ is given and $u$ is small enough then
\be{eq:non-degenerate}
I_z(\Pi^u)+\Phi(\Pi^{u}) \le z-u+u\log (u/z) +2c_\ph u <z = I_z(\d_\emptyset)+\Phi(\d_\emptyset)\,,
\ee 
so that $\d_\emptyset$ is no minimiser of the free energy.\EndProof

\bigskip
Finally we show that the energy density $\Phi$ is the infinite volume limit of 
the mean energy per volume. \bigskip

\Proof[ of Proposition \ref{prop:Phi}]
We begin with the case of periodic boundary conditions. For every $P\in\cP_\Th$, we have $ v_n\inv \int H_{n,\per} \,dP=\int \Phi(R_n)\,dP=\Phi(PR_n)$. It is easy to see that
$PR_n\to P$, cf. Remark 2.4 of \cite{GiiZess}. Since $\Phi$  is continuous,
it follows that $\Phi(PR_n)\to\Phi(P)$. 

Next we consider the case of configurational boundary conditions and suppose 
that $P$ is tempered. Applying \eqref{eq:Palm} we obtain for each $n$
\begin{eqnarray}\label{eq:P(H_n)}
\int P(d\o)\, H_{n,\o}(\o)
 &=&  
\int P(d\o)\sum_{\tau\in\sD(\o)}\ph(\tau-c(\tau))\,\1_{\{B(\tau-c(\tau))\cap(\L_n-c(\tau))\ne\emptyset\}}
\nonumber\\
&=&\int \mu_P(d\tau)\ \ph(\tau)\int dx\ \1_{\{B(\tau)\cap(\L_n-x)\ne\emptyset\}}
\\
&=&\int \mu_P(d\tau)\ \ph(\tau) \,\big|\L_n ^{\r(\tau)}\big|, \nonumber
\end{eqnarray}
where $\L_n ^{\r(\tau)}$ is the $\r(\tau)$-neigbourhood of $\L_n$. Now, 
for each $\tau$ we have 
\[
\big|\L_n ^{\r(\tau)}\big|/v_n= 1+ 4\r(\tau)/\sqrt{v_n}
+ \pi\r(\tau)^2/v_n\to 1 \text{ as }n\ti.
\]
In view of \eqref{eq:H} and \eqref{eq:ph_lb}, we can apply the dominated convergence theorem to conclude that
\[
\Phi(P) = \lim_{n\ti}  v_n\inv\int H_{n,\o}(\o)\,P(d\o),
\]
as desired.\EndProof

\subsection{The variational principle: first part}\label{subsec:var_princ1}

In this section we shall prove that each minimiser of the free energy is a 
Delaunay-Gibbs measure. We start with an auxiliary result
on the `range of influence' of the boundary condition on the events within a bounded set.
Let $\D\subset\RR^2$ be a bounded Borel set and $\o\in\O$. Writing $B(x,r)$ for the open disc in $\RR^2$ with center $x$ and radius $r$, we define
\[
\cR_\D(\o) =\big\{r>0: \ \#\,\o_{B(x,r)\setminus\D}\ge1 \text{ for all }x\in\RR^2 \text{ s.t. }B(x,r)\cap\D\ne\emptyset\big\}
\]
and $\br_\D(\o) = \inf \cR_\D(\o)$, where $\inf\emptyset :=\infty$. Let $\D^{2r}=\bigcup_{x\in\D}B(x,2r)$ be the open $2r$-neigbourhood of $\D$. We then observe the following.

\begin{lem} \label{rangefct}
(a) For all $r>0$, $\{\br_\D<r\}\in \cF_{\D^{2r}\setminus\D}$.
In particular, $\br_\D$ is $\cF_{\D^c}$-measurable.

(b) For all $P\in\cP_\Th$  we have $P(\br_\D=\infty)= P(\{\emptyset\})$.
\end{lem}
\Proof (a) Let $\tilde \cR_\D(\o)$ be defined as $\cR_\D(\o)$, except that the discs
are required to have rational centres $x\in\QQ^2$. Then 
$\cR_\D(\o)\subset\tilde \cR_\D(\o)$. Moreover, if $r<r'$ then every open $r'$-disc intersecting $\D$ contains an $r$-disc with rational center and intersecting $\D$, so that
$r'\in\cR_\D(\o)$ when $r\in \tilde \cR_\D(\o)$. This shows that
\[
\{\br_\D<r\}=\bigcup_{s\in\QQ: \; 0<s<r}\ \bigcap_{x\in\QQ^2: \; B(x,s)\cap\D\ne\emptyset}
\{N_{ B(x,s)\setminus\D}\ge1\}\,,
\]
and the last set certainly belongs to $\cF_{\D^{2r}\setminus\D}$.

(b) Since $\{\emptyset\}=\bigcap_{r\in\NN}\{N_{B(0,r)}=0\}\subset \{\br_\D=\infty\}$,
it is sufficient to prove that $\{\br_\D=\infty, N_{B(0,r)}\ge1\}$ has measure zero
for all $r>0$ and $P\in\cP_\Th$. However, if $\{\br_\D=\infty\}$ occurs then 
there exists a cone $C$
with apex at some point in the closure of $\D$ and an axis in one of finitely many
prescribed directions such that $N_C=0$, while the Poincar\'e recurrence theorem
implies that $N_C=\infty$ for each such $C$ almost surely on $\{ N_{B(0,r)}\ge1\}$.
This contradiction gives the desired result.\EndProof

\bigskip
As an immediate consequence we obtain that each nondegenerate 
stationary point process
is concentrated on the set $\O^*$ of admissible configurations.

\begin{cor} \label{temperedness} The set $\O^*$ of admissible configurations
is measurable (in fact, shift invariant and tail measurable), and
$P(\O^*)=1$ for all $P\in\cP_{\Th}$ with $P(\{\emptyset\})=0$.
\end{cor}
\Proof This is immediate from Lemma \ref{rangefct} because
$\O^* = \bigcap_{n\ge1}\{\br_{\L_n}<\infty\}$.\EndProof

\bigskip
Next we state a consequence of Proposition \ref{prop:pressure}.

\begin{cor}\label{cor:rel-entropy-wrto-Gibbs}
For every $P\in\sM_\Th(z,\ph)$, we have
\[
\lim_{n\ti} v_n\inv I_n(P;G_{n,z,\per}) = 0\;.
\]
\end{cor}
\Proof
By the definition of relative entropy, 
\[ 
I_n(P;G_{n,z,\per})= I_n(P;\Pi^z)+\int H_{n,\per} \,dP + \log Z_{n,z,\per}\,.
\] 
Together with Propositions \ref{prop:Phi} and \ref{prop:pressure}, this gives the result.\EndProof

\bigskip
We are now ready to show that the minimisers of $I_z+\Phi$ are Gibbsian.

\bigskip
\Proof[ of Theorem \ref{varprinc}, first part] 
We follow the well-known scheme of Preston 
\cite{PrestonLNM} (in the variant used in \cite{GiiJSP}, Section 7).
Let $P\in\sM_\Th(z,\ph)$, $f$ be a bounded local function, $\D$ a bounded Borel set, and 
\[
f_\D(\o) = \int f(\zeta)\,G_{\D,z,\o}(d\zeta)\,, \quad\o\in \O.
\]
We need to show that $\int f\,dP = \int f_\D\,dP$. 
Let $\br_\D$ be the range function defined above, and for each $r>0$ let 
$\1_{\D,r}= \1\{\br_\D< r\}$ and $\D^{2r}$ be the $2r$-neigbourhood of $\D$. 
By Lemma \ref{rangefct}\,(a), $\1_{\D,r}$ is measurable with respect to $\cF_{\D^{2r}\setminus\D}$. Moreover, if $\br_\D(\o)< r$ then 
$H_{\D,\o}(\zeta) = H_{\D,\o_{\D^{2r}\setminus\D}}(\zeta)$ for all
$\zeta\in\O_\D$. So, if $r$ is so large that $f$ is $\cF_{\D^{2r}}$-measurable, we can conclude that  $\1_{\D,r}\,f_\D$ is  $\cF_{\D^{2r}\setminus\D}$-measurable.

Now we apply  Corollary \ref{cor:rel-entropy-wrto-Gibbs}, which states that  $\lim_{n\ti} v_n\inv I_{\L_n}(P;G_{n,z,\per}) = 0$. 
By shift invariance, this implies that $P_{\L}\ll G_{\L,z,\per}$
with a density $g_\L$ for each sufficiently large square $\L$. 
In particular, for any $\D'\subset\L$ we have $P_{\D'}\ll (G_{\L,z,\per})_{\D'}$ wih density 
$g_{\L,\D'}(\o)=\int G_{\L\setminus\D',z,\o\cap\D'}(d\zeta)\, g_\L(\zeta)$. 
Corollary \ref{cor:rel-entropy-wrto-Gibbs} implies further that for each $\d>0$
there exists a square $\L$ and a Borel set $\D'$ with $\D^{2r}\subset\D'\subset\L$
such that $\int |g_{\L,\D'}-g_{\L,\D'\setminus\D}|\,dG_{\L,z,\per} <\d$; cf. Lemma 7.5 of
\cite{GiiJSP}. Now we consider the difference
\[
\int \1_{\D,r}\,(f-f_\D)\,dP=
\int  \1_{\D,r}\big(g_{\L,\D'}\,f-g_{\L,\D'\setminus\D}\,f_\D\big)\,dG_{\L,z,\per}\,.
\] 
Since $G_{\L,z,\per} = \int G_{\L,z,\per}(d\o)\,G_{\D,z,\o}$ and $ \1_{\D,r}\,
g_{\L,\D'\setminus\D}$ is $\cF_{\L\setminus\D}$-measurable, we can conclude that
\[
\int  \1_{\D,r}\,g_{\L,\D'\setminus\D}\,f_\D\,dG_{\L,z,\per}
=\int  \1_{\D,r}\,g_{\L,\D'\setminus\D}\,f\,dG_{\L,z,\per}\,.
\]
By the choice of $\L$ and $\D'$, we can replace the density $g_{\L,\D'\setminus\D}$
in the last expression by $g_{\L,\D'}$ making an error of at most $\d$. We thus find that
$\int \1_{\D,r}\,(f-f_\D)\,dP<\d$. Letting $\d\to0$ and $r\ti$, we finally obtain by the dominated convergence theorem that
\[
\int_{\{\br_\D<\infty\}}\,(f-f_\D)\,dP = 0\,.
\]
This completes the proof because $P(\br_\D=\infty)= P(\{\emptyset\})=0$
by Lemma \ref{rangefct}(b) and  Proposition~\ref{prop:non-degenerate}.\EndProof

\begin{rem}\label{rem:hc1} \rm
Here are some comments on the necessary modifications in the hard-core setup
of Remark \ref{rem:hc}. In analogy to Proposition \ref{prop:pressure}, one needs that
\[
\lim_{n\ti} v_n^{-1}\log  \int e^{-H_{n,\per}-H_{n,\per}^\hc} \,d\Pi^z_n
=- \min_{P\in\cP_\Th}\big [I_z(P)+\Phi(P)+\Phi^\hc(P)\big]\,.
\]
This follows directly from Propositions 4.1 and 5.4
of \cite{GiiPTRF} because $\Phi$ is continuous. Corollary 
\ref{cor:rel-entropy-wrto-Gibbs} therefore still holds for the periodic Gibbs
distributions with additional hard-core pair interaction. One also needs to modify
the proof of Proposition \ref{prop:non-degenerate}, in that the Poisson 
processes $\Pi^u$ should be replaced by the Gibbs measure $P^u$ 
with activity $u$ and pure hard-core interaction. $P^u$ is defined as 
the limit of the Gibbs distributions $G_{n,u,\per}^\hc$ for the periodic hard-core Hamiltonians $H_{n,\per}^\hc$.
By Proposition 7.4 of \cite{GiiPTRF}, $P^u$ exists and satisfies
\[
I_u(P^u)=-\lim_{n\ti} v_n^{-1}\log  \int e^{-H_{n,\per}^\hc} \,d\Pi^u_n
\le-\lim_{n\ti} v_n^{-1}\log \Pi^u_n(\{\emptyset\})= u\,.
\]
Together with \eqref{eq:I_z} we find that
\[
I_z(P^u)+\Phi(P^u)\le z+z(P^u)\, \big[ \log (u/z)+2c_\ph \big]\,,
\]
which is strictly less than $z$ when $u$ is small enough. 
Since $\Phi^\hc(P^u)=\Phi^\hc(\d_\emptyset)=0$, it follows that
the minimisers of $I_z+\Phi+\Phi^\hc$ are non-degenerate.
No further changes are required for the proof of the first part of the variational principle.
\end{rem}

\subsection{Boundary estimates}

We now work towards a proof of the reverse part of the variational principle.
In this section, we control the boundary effects that determine the difference of
$H_{n,\per}$ and $H_{n,\o}$. The resulting estimates will be crucial for the proof of Proposition~\ref{lowerbound}. 
For every $\o\in\O$ and every Borel set $\D$ let
$$ 
\sS_\D(\o)=\big\{\tau\in \sD(\o):  B(\tau)\cap \D \ne \emptyset \text{ and }  B(\tau)\setminus \D\ne \emptyset \big\}
$$
be the set of all triangles $\tau\in\sD(\o)$ for which $B(\tau)$ crosses the boundary 
of $\D$. We start with a lemma that controls the influence on $\sS_\D$ when two
configurations are pasted together.

\begin{lem}\label{frontiere}
Let $\Delta$ be a (not necessarily bounded) Borel set in $\RR^2$,
$\zeta\in\O^*\cup\{\emptyset\}$ a configuraton with $\zeta_{\partial\D}=\emptyset$,
and $\o\in\O$.
Then for each
 $\tau \in \sS_\D(\zeta_\D \cup \o_{\D^c})$  and each $x\in \tau_\D$ there exists some $\tau' \in \sS_\D(\zeta)$ with $x \in \tau'$. 
\end{lem}
\Proof
Let $\D$, $\zeta$ and $\o$ be given. If $\zeta$ is empty, there exists no 
$x\in \tau_\D\subset\zeta_\D$, so that the statement is trivially true.  
So let $\zeta\in\O^*$ and suppose there exists some
$\tau\in \sS_\D(\zeta_\Delta \cup \o_{\Delta^c})$ with $\tau_\D\ne\emptyset$.
Let $x\in\tau_\D$. Since $\zeta_{\partial\D}=\emptyset$,
$x$ does in fact belong to the interior of $\D$. This implies that 
$B(\tau')\cap \D\ne \emptyset$ for each $\tau'\in\sD(\zeta)$ containing $x$.
Therefore we only need to show that $B(\tau')\setminus \D\ne \emptyset$
for at least one such $\tau'$. Suppose the contrary. Then $B(\tau')\subset \D$
whenever $x\in\tau'\in\sD(\zeta)$. This means that the Delaunay triangles
containing $x$ are completely determined by $\zeta_\D$. This gives the contradiction
\[
\emptyset\ne\{\tau\in \sS_\D(\zeta_\D \cup \o_{\D^c}):\tau\ni x\} 
= \{\tau'\in \sS_\D(\zeta):\tau'\ni x\} =\emptyset\,,
\]
and the proof is complete.\EndProof

\bigskip
The following proposition is the fundamental boundary estimate. 
It bounds the difference of Hamiltonians with periodic and configurational 
boundary conditions in terms of $S_n:=  \#\sS_{\L_n}$.

\begin{prop}\label{boundary}
 There exists a universal constant $\g<\infty$ such that 
 \[  |H_{n,\per}(\zeta)-H_{n,\o}(\zeta)|  \le \g c_\ph\,\big(S_n(\o)+S_n(\zeta)\big) \]
 for all $n\ge 1$ and all $\zeta,\o \in \O^*\cup\{\emptyset\}$ with $\zeta_{\partial\L_n}
 = \o_{\partial\L_n}=\emptyset$. 
\end{prop}

\Proof
Let $n,\zeta,\o$ be fixed and 
\[
\sA= \big\{\tau\in\sD(\zeta_{\L_n}\cup\o_{\L_n^c}):
 \, B(\tau)\cap\L_n\ne\emptyset\big\} \,,\quad
\sB= \big\{\tau\in\sD(\zeta_{n,\per}): \, c(\tau)\in\L_n\big\}\,.
\]
In view of (\ref{eq:cH_o}) and (\ref{eq:H_per}) we have
$ H_{n,\o}(\zeta)= \sum_{\tau\in\sA}\ph(\tau)$ and 
$ H_{n,\per}(\zeta)=\sum_{\tau\in\sB}\ph(\tau)$.
Since $\ph$ is bounded by $c_\ph$, we only need to estimate
the cardinalities of $\sA\setminus\sB$ and $\sB\setminus \sA$. 
We note that
$\sA\setminus \sB \subset \sS_{\L_n}(\zeta_{\L_n}\cup\o_{\L_n^c})$ and
$\sB\setminus\sA \subset\sS_{\L_n}(\zeta_{n,\per})$.
So we can apply Lemma \ref{frontiere} to both $\D=\L_n$ and $\D=\L_n^c$ 
to obtain that 
the set of points belonging to a triangle in $\sA\setminus \sB$ is contained in the set of points belonging to a triangle of $\sS_{\L_n}(\zeta)\cup \sS_{\L_n}(\o)$. Hence,
$\#\big(\bigcup_{\tau\in\sA\setminus \sB}\tau\big) \le 3(S_n(\zeta) + S_n(\o))$.
By Lemma \ref{Euler}, it follows that $\#(A\setminus B) \le 6 (S_n(\o) + S_n(\zeta))$.

To estimate the cardinality of $\sB \setminus\sA$ we may assume that $\zeta_n\ne\emptyset$. The periodic continuation $\zeta_{n,\per}$ then contains a lattice, and this
implies that every triangle of $\sD(\zeta_{n,\per})$ has a circumscribed disc of diameter
at most $\sqrt{2v_n}$. Hence, each $\tau\in \sB\setminus \sA$ is contained in 
$\L_{5n+2}$, the union of $5^2$ translates of $\L_n$ (up to their boundaries). 
Applying Lemma~\ref{frontiere} to each of these translates we conclude that the 
number of points that belong to a triangle of $\sB\setminus \sA$ is bounded by 
$3\cdot5^2\, S_n(\zeta)$. Using Lemma \ref{Euler} again we find that
$\#(B\setminus A) \le 150 (S_n(\o) + S_n(\zeta))$,
and the result follows with $\g=156$.\EndProof

\bigskip
The following immediate corollary will be needed in the proof of Theorem \ref{varprinc}.

\begin{cor}\label{vide}
 There exists a constant $C<\infty$ such that
$ | H_{n,\o}(\emptyset)|  \le C\, S_n(\o)$  for all $n\ge 1$ and $\o\in\O^*$
with $\o_{\partial\L_n}=\emptyset$.
\end{cor}

The next proposition exhibits the fundamental role of the temperedness 
condition~\eqref{eq:H} combined with stationarity for controlling the boundary effects.

\begin{prop}\label{convergence}
For every $P\in\cP_\Th^\tmp$, $v_n\inv S_n\to 0$ in $L^1(P)$ and $P$-almost surely.
\end{prop}

\Proof 
For each $i\in\RR^2$ 
we consider the shifted unit square $C(i)=\L_0+i$ and define the random variable
\[
Z_i=\#\big\{\tau\in\sD(\cdot):\, B(\tau)\cap C(i)\ne \emptyset\big\}.
\]
Then 
\be{eq:S-bound} S_n \le \sum_{i\in I_n\setminus I_{n-1}} Z_i\,,\ee
where $I_n=\L_n\cap(\ZZ^2+(\frac12,\frac12))$.  Note that $\#\, I_n=v_n$.
As in \eqref{eq:P(H_n)} we have
\[
\int Z_0\,dP =  
\int \mu_P(d\tau)\ \big|\L_0^{\r(\tau)}\big|
=\int \mu_P(d\tau)\ \big(1+ 4\r(\tau)+ \pi\r(\tau)^2\big)\,.
\]
The last term is finite by the temperedness of $P$. 
So, each $Z_i$ is $P$-integrable.
Since $Z_i=Z_0\circ\th_i$, the two-dimensional ergodic theorem
implies that $v_n\inv\sum_{i\in I_n} Z_i$ converges to a finite limit $\bar Z$,
both $P$-almost surely and in $L^1(P)$. This implies that
$v_n\inv \sum_{i\in I_n\setminus I_{n-1}} Z_i$ tends to zero $P$-almost surely and in $L^1(P)$. The result thus follows from \eqref{eq:S-bound}.\EndProof

\subsection{Temperedness and block average approximation}\label{subsec:block_av}

Our first result in this subsection is a sufficient condition for temperedness in terms
of vacuum probabilities. For $P\in\cP_\Th$ let
\be{eq:V_k}
V_k(P) = \text{ess\,sup}\; P(N_{\L_k}=0|\cF_{\L_k^c})
\ee
be the essential supremum of the conditional probability that $\L_k$ contains
no particle given the configuration outside.

\begin{prop}\label{tempered}
Every $P\in\cP_\Th$ satisfying
\begin{equation}\label{empty}
\sum_{k\ge0}v_k \,V_k(P)\,z(P) <\infty
 \end{equation}
 is tempered.
\end{prop}
\Proof 
We can assume that $P\ne\d_\emptyset$ because otherwise the result is trivial.
For each  $k\ge 1$  we consider the shifted squares 
$\L_k(i)=\L_k+(2k+1)i$, $i\in\ZZ^2$, as well as the event
\[
 A_k=\big\{N_{\L_k(i)} \ne0 \text{ for all $i\in\ZZ^2$ with } \|i\|_\infty\le1 \big\}.
 \]
Since $P\ne\d_\emptyset$, it is clear that $P(A_k)\to 1$ as $k\ti$. Thus we can write
\[\begin{split}
\int \#&\big(\tau\in\sD(\o): 0\in B(\tau)\big) \,P(d\o) \\
&\le \ \sum_{k\ge 1} \int_{A_k \setminus A_{k-1}} 
\#\big(\tau\in\sD(\o): 0\in B(\tau)\big) P(d\o)
\end{split}\]
with the convention $A_0=\emptyset$. Now, if $\o\in A_k$ then each circumscribed disc
containing $0$ of a triangle $\tau\in \sD(\o)$ has a diameter not larger than  
$2\sqrt{2v_k}$, so that each such $\t$ in fact belongs to $\sD(\o_{\L_{7k+3}})$.
Lemma \ref{Euler} thus shows that the number of such $\t$ is at most
$2N_{7k+3}(\o)$.
The last sum is therefore not larger than
\[
\sum_{k\ge 1} \int_{A_{k-1}^c} 2N_{7k+3}\, dP
 \le  2 \sum_{k\ge 1} \sum_{ i\in\ZZ^2: \|i\|_\infty\le1} \int \1_{\{N_{\L_{k-1}(i)}=0\}} 
N_{7k+3}\, dP \,.
\]
In view of the stationarity of $P$, the last integral is bounded by
$V_{k-1}(P)\, v_{7k+3}\,z(P)=7^2\,v_k\, V_{k-1}(P)\, z(P)$. So we arrive at the estimate
\[
\int \#\big(\tau\in\sD(\o): 0\in B(\tau)\big) \,P(d\o)\le
2\cdot 7^2\,3^2\sum_{k\ge 1}v_k\, V_{k-1}(P)\, z(P) \,.
\]
Together with \eqref{eq:tempered2} and assumption \eqref{empty}, this  implies the temperedness of $P$ because $ v_k\sim v_{k-1}$ as $k\ti$.\EndProof

\bigskip
The second result concerns the approximation of stationary measures in terms
of tempered ergodic measures. This approximation uses the block average construction
first introduced by Parthasarathy for proving that the ergodic measures are dense in 
$\cP_\Th$; cf. \cite[Theorem (14.12)]{GiiGibbs}, for example.

\begin{prop}\label{ergodic_approx} Let $z>0$ and $Q\in\cP_\Th$ be such that $I_z(Q)+\Phi(Q)<\infty$. Then
for each $\e>0$ there exists some tempered $\Th$-ergodic $\hat Q\in\cP_\Th$ such that $I_z(\hat Q)<I_z(Q)+\e$ and $\Phi(\hat Q)<\Phi(Q)+\e$.
\end{prop}
\Proof
Let $Q\in\cP_\Th$ be given.
We can assume that $Q\ne\d_\emptyset$ because otherwise we can choose
$\hat Q=Q$.
For $n\ge1$ let  $Q^{\text{iid}}_n$ 
denote the probability measure on $\O$ relative
to which the particle configurations in the blocks $\L_n+
(2n+1)i$, $i\in \ZZ^2$, are independent with identical distribution $Q_n=Q\circ \pr_{\L_n}\inv$. 
(In particular, this means that the boundaries of these blocks contain no particles.)
Consider the spatial average
\[
Q^{\text{iid-av}}_n=v_n\inv\int_{\L_n}Q^{\text{iid}}_n\circ\th_x^{-1}dx.
\]
It is clear that $Q^{\text{iid-av}}_n\in\cP_\Th$. It is also well-known that 
$Q^{\text{iid-av}}_n$ is $\Th$-ergodic; cf.~\cite[Theorem (14.12)]{GiiGibbs}, for example.
By an analogue of \cite[Proposition (16.34)]{GiiGibbs} or \cite[Lemma 5.5]{GiiZess}, 
its entropy density satisfies
\be{eq:block_entropy}
I_z(Q_n^{\text{iid-av}})=v_n\inv I(Q;\Pi^z_n).
\ee
So,  $I_z(Q_n^{\text{iid-av}})< I_z(Q)+\e$ when $n$ is large enough.
Moreover, the same argument as in \cite[Lemma 5.7]{GiiZess} shows that
$Q_n^{\text{iid-av}}\to Q$ in the topology $\Tau_\cL$. 
By Proposition \ref{muP_cont}, $\Phi$ is continuous. We thus conclude
that $\Phi(Q_n^{\text{iid-av}})\to \Phi(Q)$,  whence $\Phi(Q_n^{\text{iid-av}})<\Phi(P)+\e$
for large~$n$.
It remains to prove that each $Q_n^{\text{iid-av}}$ is tempered.  
Let $n$ be fixed and $k\ge\ell(2n{+}1)$ for some $\ell\ge1$. 
We claim that $V_k(Q_n^{\text{iid-av}})\le q^{v_{\ell-1}}$ with
$q=Q(N_n=0)$. 
Indeed, for each $x\in\L_n$ we have $\L_k+x\supset\L_{n+(\ell-1)(2n+1)}$, and 
the latter set consists of
$v_{\ell-1}=(2\ell-1)^2$ distinct blocks as above.
Letting $g$ be any nonnegative $\cF_{\L_k^c}$-measurable function
and using  the independence of block configurations, we thus 
conclude that
\[\begin{split}
&\int \1_{\{N_k=0\}}g\,d Q_n^{\text{iid-av}}
=v_n\inv\int_{\L_n}dx\int dQ_n^{\text{iid}} \, \1_{\{N_{\L_k+x}=0\}}\,g\circ\th_x\\
&\le v_n\inv\int_{\L_n}dx\int dQ_n^{\text{iid}} \, \1_{\{N_{n+(\ell-1)(2n+1)}=0\}}\,g\circ\th_x
=q^{v_{\ell-1}} \int g\,d Q_n^{\text{iid-av}}\,,
\end{split}\]
which proves the claim. Now, we have $q<1$ because $Q\ne\d_\emptyset$. It follows that
\[
\sum_{k> 2n}v_k\,V_k(Q_n^{\text{iid-av}})
\le \sum_{\ell\ge1}q^{v_{\ell-1}}\sum_{\ell(2n+1)\le k<(\ell+1)(2n+1)}v_k
\le C_n  \sum_{\ell\ge0}v_\ell \,q^{v_\ell}<\infty
\]
for some constant $C_n<\infty$. Together with Proposition~\ref{tempered},
this gives the temperedness of $Q_n^{\text{iid-av}}$.\EndProof

\subsection{The variational principle: second part}\label{subsec:var_princ2}

In this section we will complete the proof of the variational principle.
The essential ingredient is the following counterpart of Proposition \ref{prop:pressure}
involving configurational instead of periodic boundary conditions. We only state 
the lower bound we need.

\begin{prop}\label{lowerbound}
For every $P\in\cP_\Th^\tmp$ and $P$-almost every $\o$ we have
\[
\liminf_{n\ti}v_n\inv\log Z_{n,z,\o} \ge p(z,\ph)\,.
\]
\end{prop}
\Proof 
By \eqref{eq:pressure} and Lemma \ref{ergodic_approx}, it is sufficient to show that
\[
\liminf_{n\ti}v_n\inv\log Z_{n,z,\o} \ge -I_z(Q)-\Phi(Q)
\] 
for every ergodic $Q\in \cP_\Th^\tmp$.
We can assume without loss that the right-hand side is finite. Now, since $I_z(Q)$ is finite, $Q$ is locally absolutely continuous with repect to $\Pi^z$. 
We fix some $\e>0$, let $f_n=dQ_n/d\Pi^z_n$, and consider for every $\o\in\O^*$ the set
\[
A_{n,\o}=\big\{|H_{n,\o}-H_{n,\per}|/v_n\le \e,\, \Phi(R_n)<\Phi(Q)+\e,\, v_n\inv\log f_n< I_z(Q)+\e\big\}\,.
\]
Then for sufficiently large $n$ we have
\begin{eqnarray*}
Z_{n,z,\o} & \ge & \int_{A_{n,\o}} e^{-H_{n,\o}} \,f_n\inv  \,d Q \\
& \ge & \int_{A_{n,\o}} e^{-H_{n,\per}}e^{-v_n\e} \,f_n\inv  \,d Q \\
& \ge &  e^{-v_n\,[I_z(Q) +\Phi(Q)+3\e]}\, Q(A_{n,\o})\,. 
\end{eqnarray*}
It is therefore sufficient to show that for $P$-almost $\o$, $Q(A_{n,\o})\to 1$ as $n\ti$.
By the ergodic theorem, $\Phi(R_n)$ converges to $\Phi(Q)$ in $Q$-probability; cf. Remark 2.4 of \cite{GiiZess}. By McMillan's theorem \cite{Fritz, NgZess}, $Q(v_n\inv\log f_n < I_z(Q)+\e)\to 1$ when $n$ tends to infinity. Moreover, Propositions \ref{boundary} and \ref{convergence} imply that, for $P$-almost all $\o$, $|H_{n,\o}-H_{n,\per}|/v_n$ converges to $0$ in $L^1(Q)$. This gives the result.\EndProof

\bigskip
We can now show that every tempered stationary Gibbs measure minimises the free energy density.

\bigskip
\Proof[ of Theorem \ref{varprinc}, second part]
We follow the argument of \cite{GiiJSP}, Proposition~7.7. Let $P\in\sG^\tmp_\Th(\ph,z)$.
On each $\cF_{\L_n}$, $P$ is absolutely continuous w.r. to $\Pi^z$ with density
\[
g_n(\zeta) = \int P(d\o)\, \frac{d G_{n,z,\o}}{d\Pi^z_n}(\zeta) =
 \int P(d\o)\, e^{-H_{n,\o}(\zeta)}/ Z_{n,z,\o} \,.
\]
Using Jensen's inequality and the Gibbs property of $P$ we thus find that
\[\begin{split}
I_n(P;\Pi^z) &= \int  g_n \log g_n \, d\Pi^z\\
&\le \int \Pi^z(d\zeta)\int P(d\o) \frac{d G_{n,z,\o}}{d\Pi^z_n}(\zeta) 
\big[-H_{n,\o}(\zeta)-\log Z_{n,z,\o}\big]\\
&= - \int P(d\o)\, H_{n,\o}(\o) - \int P(d\o)\,\log Z_{n,z,\o}\,.
\end{split}\]
Next we divide by $v_n$ and let $n\ti$.
We know from Proposition \ref{prop:Phi} that $v_n\inv \int P(d\o)\, H_{n,\o}(\o)\to \Phi(P)$.
On the other hand,  Corollary \ref{vide} implies that $v_n\inv\,\log Z_{n,z,\o}\ge -z- Cv_n\inv S_n(\o)$. Using Propositions \ref{convergence} and \ref{lowerbound} 
together with Fatou's Lemma, we thus find that
\[\liminf_{n\ti} v_n\inv \int P(d\o)\,\log Z_{n,z,\o} \ge p(z,\ph).\]
Therefore
$I_z(P) \le -\Phi(P) + \min[I_z+\Phi]$,  as required.\EndProof

\begin{rem}\label{rem:hc2} \rm
In the hard-core setting of Remark \ref{rem:hc}, a slight refinement of Proposition \ref{ergodic_approx} is needed. Namely, under the additional assumption that $\Phi^\hc(Q)=0$ one needs to achieve that also $\Phi^\hc(\hat Q)=0$.
To this end we fix an integer $k>r_0/2$ and define $Q^{\text{iid}}_n$ 
in such a way that the particle configurations in the blocks $\L_n+
(2n+1)i$, $i\in \ZZ^2$, are independent with identical distribution 
$Q_{n-k}$, rather than $Q_n$. This means that the blocks are separated by corridors
of width $2k>r_0$ that contain no particles. It follows that $\Phi^\hc(Q_n^{\text{iid-av}})=0$, 
and it is still true that $\limsup_{n\ti}I_z(Q_n^{\text{iid-av}})\le I_z(Q)$; cf. 
\cite[Lemma 5.1]{GiiPTRF}. We thus obtain the refined Proposition \ref{ergodic_approx}
as before. 

A similar refinement is required in the proof of Proposition \ref{lowerbound}.
One can assume that $\Phi^\hc(P)=0$ and $\Phi^\hc(Q)=0$, and in the definition
of $A_{n,\o}$ one should introduce an empty corridor at the inner boundary
of $\L_n$ to ensure that $H^\hc_{n,\o}=H^\hc_{n,\per}=0$ on $A_{n,\o}$, see \cite[Proposition 5.4]{GiiPTRF} for details.
In the proof of the second part of Theorem~\ref{varprinc}, one then only needs to note
that $\Phi^\hc(P)=0$  when $P$ is a Gibbs measure $P$ for the combined
triangle and hard-core pair interaction.
The proof of Theorem \ref{Gtempered}
carries over to the hard-core case without any changes.
\end{rem}

\subsection{Temperedness of Gibbs measures}\label{Gibbstempered}

Here we establish Theorem \ref{Gtempered}. By Proposition \ref{tempered}
it is sufficient to show the following.
\begin{prop}\label{temperedb}
Let $\ph$ be bounded and eventually increasing, $z>0$, and
$P$ be any stationary Gibbs measure for $\ph$ and $z$. Then there exists a constant 
$C>0$ such that
\begin{equation}\label{emptyb}
 P(N_k=0 | \cF_{\L_k^c}) \le {C}\,{v_k^{-2}}
 \end{equation}
for all $k\ge 1$.
\end{prop}

\begin{figure}
\begin{center}\small
\psfrag{x}{\hspace{-1.5mm}$x$}
\epsfig{file=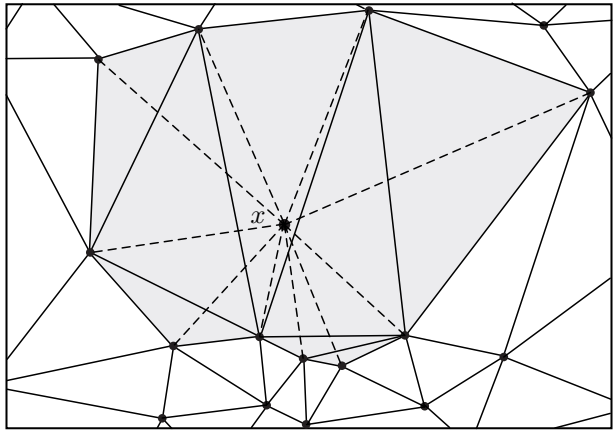, width=8cm}
\caption{$\sD(\o)$ (solid lines) and $\sD(\o\cup\{x\})
\setminus\sD(\o)$ (dashed lines).
The difference region $\D_x(\o)$ is shaded in grey. \label{DiffGraph}}
\end{center}
\end{figure}

To prove this we need an auxiliary result which states that the radii of all circumcircles in the Delaunay tessellation must decrease when a point is added to the configuration.
Specifically, let $\o\in \O^*$ and $x\in \RR^2\setminus\o$ be such that $\o\cup\{x\}$ is in general circular position and $x$ is not collinear with two points of $\o$. 
We consider the sets
\[
\sC_x(\o):=\sD(\o) \setminus \sD(\o\cup\{x\})=\big\{\tau\in\sD(\o): B(\tau)\ni x\big\}
\]
and 
\[ 
\sC^+_x(\o):=\sD(\o\cup\{x\}) \setminus \sD(\o)=
\big\{\tau\in\sD(\o\cup\{x\}):\tau\ni x\big\} . 
\]
If $\langle\tau\rangle$ denotes the convex hull of a triangle $\tau$,
\be{eq:diff_region}
\D_x(\o):=\bigcup_{\tau\in\sC_x(\o)}\langle\tau\rangle=
\bigcup_{\tau\in\sC_x^+(\o)}\langle\tau\rangle
\ee
is the region on which the triangulations $\sD(\o)$ and $\sD(\o\cup\{x\})$ differ;
see Fig.~\ref{DiffGraph}. 
Up to the point $x$, the interior $\D_x^o(\o)$ of $\D_x(\o)$ is covered by the discs $B(\tau)$ with $\tau\in\sC_x^+(\o)$, which by definition are free of particles.
Consequently, $\D_x^o(\o)$ contains no particle of $\o$, so that the vertices of each
$\tau\in\sC_x(\o)$ belong to the boundary $\partial\D_x(\o)$.
Next, Lemma \ref{Euler} shows that
\begin{equation}\label{card}
\# \sC^+_x(\o)= \# \sC_x(\o)+2\,,
\end{equation}
and for every $\L\ni x$ we have 
\begin{equation}\label{energielocale}
H_{\L,\o}(\o\cup\{x\})- H_{\L,\o}(\o) = 
\sum_{\tau \in \sC^+_x(\o)} \ph(\tau) - \sum_{\tau \in \sC_x(\o)} \ph(\tau).
\end{equation}
Here is the monotonicity result announced above.

\begin{lem}\label{ordre} Under the conditions above, there exist
a subset $\sC_x'(\o)\subset \sC_x(\o)$ with 
$\#\big(\sC_x(\o)\setminus \sC_x'(\o)\big)\le 4$
and an injection $I$ from $\sC_x'(\o)$ to $\sC^+_x(\o)$ such that 
\begin{equation}\label{inegalite}
\r(I(\tau)) \le \r(\tau)\quad\text{ for all $\tau\in \sC_x'(\o)$.}
\end{equation}
\end{lem}
We postpone the proof of this lemma until the end, coming first to its use. 

\bigskip
\Proof[ of Proposition \ref{temperedb}]
By assumption, $\ph$ is eventually increasing. So there exists some $r_\ph<\infty$
and a nondecreasing function $\psi$ such that $ \ph(\tau) =\psi(\r(\tau))$ when 
$\r(\tau) \ge r_\ph$. Combining Lemma~\ref{ordre} and Equations (\ref{card}) and (\ref{energielocale})
we thus find that
\begin{equation}\label{energie}
 H_{k,\o}(\o\cup\{x\}) \le  H_{k,\o}(\o) + 10\,c_\ph
 \end{equation}
for all $\o\in \O^*$, $k\ge1$, and Lebesgue-almost all $x\in \L_k\setminus\o$ that have
at least the distance $2r_\ph$ from all points of $\o$. Next, let $P\in\sG_\Th(z,\ph)$. 
By definition, 
\[
P(N_k=0 | \cF_{\L_k^c})(\o)=Z_{k,z,\o}\inv\, e^{-zv_k} \,e^{-H_{k,\o}(\emptyset)}
\]
for all $\o\in\O^*$. Let 
$\L_k^{(2)}=\big\{(x,y)\in\L_{k-2r_\ph}^2: |x-y|\ge 2r_\ph\big\}$.
Applying (\ref{energie}) twice (viz. to $\o_{\L_k^c}$ and $x$ as well as $\o_{\L_k^c}\cup\{x\}$ and $y$)
and recalling \eqref{eq:Z_L} we find that
\[\begin{split}
Z_{k,z,\o} 
&\ge e^{-zv_k}\, \frac{z^2}{2}  \int_{\L_k^{(2)}} e^{-H_{k,\o}(\{x\}\cup\{y\})}\,dx\, dy\\  
&\ge {z^2\,|\L_k^{(2)}|}\,e^{-zv_k}\, e^{-H_{k,\o}(\emptyset)-20\,c_\ph}/\,2\,.
\end{split}\] 
Since $|\L_k^{(2)}|\sim v_k^2$ as $k\ti$, the result follows.\EndProof

\bigskip
Finally we turn to the proof of Lemma \ref{ordre}.
\bigskip

\Proof[ of Lemma \ref{ordre}] 
Let $\tau_x$ be the unique triangle of $\sC_x(\o)$ containing $x$ in its interior,
and $\sC_x^{+\wedge}(\o)$ the set of all $\tau\in\sC_x^+(\o)$ that have an acute
or right angle at $x$. Note that $\#(\sC_x^+(\o)\setminus\sC_x^{+\wedge}(\o))\le 3$
because the angles at $x$ of all $\tau\in\sC_x^+(\o)$ add up to 360 degrees.
We will associate to each triangle $\tau\in\sC_x(\o)$ a triangle $I(\tau)\in\sC^+_x(\o)$,
except possibly when $\tau=\tau_x$ or the
candidate for $I(\tau)$ does not belong to $\sC_x^{+\wedge}(\o)$.
Our definition of $I(\tau)$ depends on the number $k=k(\tau)$ of edges $e\subset\tau$
with $\langle e\rangle\subset\partial\D_x(\o)$. 
Let $\sC^{(k)}_x(\o)$ be  the set of all $\tau\in\sC_x(\o)$ that have
$k$ such edges.
Since $\sC^{(3)}_x(\o)=\emptyset$ except when $\sC_x(\o)=\{\tau_x\}$, 
we only need to consider the three cases $k=0,1,2$.

The cases $k=1$ and $2$ are easy: For every $\tau\in\sC^{(1)}_x(\o)$ there exists a unique edge $e(\tau)$ such that $e(\tau)\cup\{x\}\in\sC_x^+(\o)$. If in fact 
$e(\tau)\cup\{x\}\in\sC_x^{+\wedge}(\o)$ we set $I(\tau)=e(\tau)\cup\{x\}$; 
otherwise we leave
$I(\tau)$ undefined. Likewise, every $\tau\in \sC^{(2)}_x(\o)$ has two edges 
$e_1(\tau)$ and $e_2(\tau)$ in $\partial\D_x(\o)$ (in clockwise order, say) 
and can be mapped to $I(\tau)=e_1(\tau)\cup\{x\}$, provided this triangle 
belongs to $\sC^{+\wedge}_x(\o)$. The resulting mapping $I$ is clearly injective. 
Moreover, $\tau$ and $I(\tau)$ have the edge $e(\tau)$ (resp.~$e_1(\tau)$) in common,
and $x\in B(\tau)$ because $\tau\in\sC_x(\o)$. Since $I(\tau)\in\sC^{+\wedge}_x(\o)$
whenever it is defined, we can conclude that $\r(I(\tau)) \le \r(\tau)$.

The case $k=0$ is more complicated because the tiles $\tau\in\sC^{(0)}_x(\o)$ are not naturally associated to a tile of $\sC^+_x(\o)$. 
To circumvent this difficulty we define an injection $\tilde I$ from 
$\sC^{(0)}_x(\o)\setminus \{\tau_x\}$ to $\sC^{(2)}_x(\o)$ such that 
$\r(\tilde I(\tau)) \le \r(\tau)$. Each triangle $\tau\in\sC^{(0)}_x(\o)$
different from $\tau_x$ can then be mapped to the triangle $I(\tau)=e_2(\tilde I(\tau))\cup\{x\}$,
provided the latter belongs to $\sC_x^{+\wedge}(\o)$; otherwise $I(\tau)$ remains
undefined. This completes the construction of $I$.
(Note that $\tau_x$ does not necessarily belong to $\sC^{(0)}_x(\o)$. However, if
it does we have no useful definition of $\tilde I(\tau_x)$.)

To construct $\tilde I$ we turn $\sC_x(\o)$ into the vertex set of a graph $G_x(\o)$
by saying that two tiles are adjacent if they share an edge. 
The set $\sC^{(2)}_x(\o)$ then coincides with the set of all leaves of $G_x(\o)$, and 
$\sC^{(0)}_x(\o)$ is the set of all triple points of $G_x(\o)$. 
\begin{figure}
\begin{center}\small
\psfrag{x}{$x$}
\psfrag{e1}{$e_1$}
\psfrag{e2}{$e_2$}
\psfrag{t}{$\tau_*$}
\psfrag{tt}{$\tau_1$}
\psfrag{Itt}{$\tau_2$}
\psfrag{z0}{$z_0$}
\psfrag{z1}{$z_1$}
\psfrag{z2}{$z_2$}
\epsfig{file=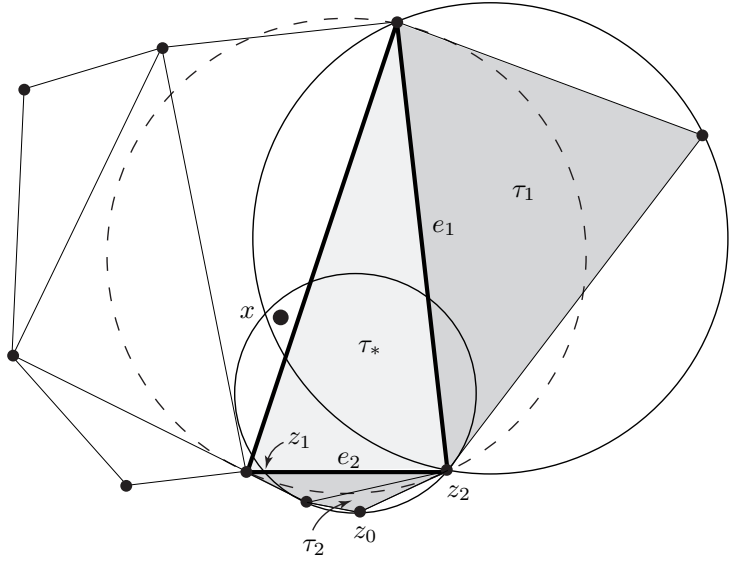,width=9.5cm}
\caption{The set $\sC_x(\o)$ for the configuration $\o$ of Fig.\,1, 
with a tile $\tau_*\in\sC_x^{(0)}(\o)$ (light grey), its circumcircle $B(\tau_*)$ (dashed), 
the associated edges $e_i$ and regions $W_i$ (dark grey), and two triangles $\tau_i\in\sC_x(\o)$ with $\tau_i\subset W_i$ with their circumcircles (solid). 
The construction in the proof gives $\tilde I(\tau_*)=\tau_2$.
\label{fig:circles}}
\end{center}
\end{figure}
Consider a fixed $\tau_*\in\sC^{(0)}_x(\o)\setminus \{\tau_x\}$. Since $\tau_*\subset \partial\D_x(\o)$,
the set $\D_x(\o)\setminus\langle\tau_*\rangle$ splits into three connected components.
Let $W_i=W_i(\tau_*,x,\o)$ be the closure of the $i$th component, $i=1,2,3$.
Any two of these sets intersect at a point of $\tau_*$, and one of them contains
$x$ because $\tau_*\ne\tau_x$. Suppose $x\in W_3$. For $i=1,2$ let 
$e_i=\tau_*\cap W_i$ be the edge of $\tau_*$ that separates $W_i$ from the 
rest of $\D_x(\o)$; see Fig.~\ref{fig:circles}.
We claim that there exists some $i=i(\tau_*)\in\{1,2\}$ such that
\be{eq:claim}
\r(\tau)\le\r(\tau_*)\quad\text{ for all $\tau\in\sC_x(\o)$ with $\tau\subset W_i$}\,.
\ee
The image $\tilde I(\tau_*)$ of $\tau_*$  can then be defined as the leaf of $G_x(\o)$
in $W_{i(\tau_*)}$ with the largest circumcircle. (The largest circumcircle condition
takes account of the fact that the path in $G_x(\o)$ from $\tau_*$ to 
$\tilde I(\tau_*)$ might contain further triple points.)

It remains to prove \eqref{eq:claim}. Since $\tau_*\ne\tau_x$, there exists 
at least one $i$ such that the triangle $\{x\}\cup e_i$ has an acute angle at $x$.
We fix such an $i$ and consider any $\tau\in\sC_x(\o)$ with $\tau\subset W_i$.
There exists at least one point $z_0\in\tau$ that is not contained in the closed
disc $\overline{B}(\tau_*)$. Since $x\in B(\tau)$ and $\langle\tau\rangle$
is covered by the tiles $\langle\tau'\rangle$ for $\tau'\in\sC^+_x(\o)$ with 
$\langle\tau'\rangle \cap \langle e_i\rangle\ne\emptyset$, we conclude that
the line segment $s$ from $z_0$ to $x$ is contained in $\overline{B}(\tau)$
and hits both the circle $\partial B(\tau_*)$ and the edge $\langle e_i\rangle$.
In particular, $B(\tau)\cap\langle e_i\rangle\ne\emptyset$. 
Since $B(\tau)$ contains no points of $\o$,
we deduce further that the circle $\partial B(\tau)$ hits the edge $\langle e_i\rangle$ in
precisely two points $z_1$ and $z_2$. By the choice of~$i$, the angle of the triangle
$\{z_1,x,z_2\}$ at $x$ is acute. Since $x\in B(\tau)$, it follows that the angle
of the triangle $\{z_1,z_0,z_2\}$ at $z_0$ is obtuse. Consequently, if we 
consider running points $y_k$ such that $y_0$ runs from
$z_0$ to the point $s\cap\partial B(\tau_*)$ and the edge $\{y_1,y_2\}$ from 
$\{z_1,z_2\}$ to $e_i$, the associated circumcircles $B(\{y_1,y_0,y_2\})$
run from $B(\tau)$ to $B(\tau_*)$,  and their radii $\r(\{y_1,y_0,y_2\})$ must grow.
This proves that $\r(\tau)\le \r(\tau_*)$,
and the proof of  \eqref{eq:claim} and the lemma is complete.\EndProof

\bigskip 
{\bf Acknowledgment}. We are grateful to Remy Drouilhet who brought
us together and drew the interest of H.-O.\,G. to the subject.


\end{document}